\catcode`\@=11
\def\undefine#1{\let#1\undefined}
\def\newsymbol#1#2#3#4#5{\let\next@\relax
 \ifnum#2=\@ne\let\next@\msafam@\else
 \ifnum#2=\tw@\let\next@\msbfam@\fi\fi
 \mathchardef#1="#3\next@#4#5}
\def\mathhexbox@#1#2#3{\relax
 \ifmmode\mathpalette{}{\m@th\mathchar"#1#2#3}%
 \else\leavevmode\hbox{$\m@th\mathchar"#1#2#3$}\fi}
\def\hexnumber@#1{\ifcase#1 0\or 1\or 2\or 3\or 4\or 5\or 6\or 7\or 8\or
 9\or A\or B\or C\or D\or E\or F\fi}

\newdimen\ex@
\ex@.2326ex
\def\varinjlim{\mathop{\vtop{\ialign{##\crcr
 \hfil\rm lim\hfil\crcr\noalign{\nointerlineskip}\rightarrowfill\crcr
 \noalign{\nointerlineskip\kern-\ex@}\crcr}}}}
\def\varprojlim{\mathop{\vtop{\ialign{##\crcr
 \hfil\rm lim\hfil\crcr\noalign{\nointerlineskip}\leftarrowfill\crcr
 \noalign{\nointerlineskip\kern-\ex@}\crcr}}}}
\def\varliminf{\mathop{\underline{\vrule height\z@ depth.2exwidth\z@
 \hbox{\rm lim}}}}

\font\tenmsa=msam10
\font\sevenmsa=msam7
\font\fivemsa=msam5
\newfam\msafam
\textfont\msafam=\tenmsa
\scriptfont\msafam=\sevenmsa
\scriptscriptfont\msafam=\fivemsa
\edef\msafam@{\hexnumber@\msafam}
\mathchardef\dabar@"0\msafam@39
\def\dashrightarrow{\mathrel{\dabar@\dabar@\mathchar"0\msafam@4B}}
\def\dashleftarrow{\mathrel{\mathchar"0\msafam@4C\dabar@\dabar@}}

\font\tenmsb=msbm10
\font\sevenmsb=msbm7
\font\fivemsb=msbm5
\newfam\msbfam
\textfont\msbfam=\tenmsb
\scriptfont\msbfam=\sevenmsb
\scriptscriptfont\msbfam=\fivemsb
\edef\msbfam@{\hexnumber@\msbfam}
\def\Bbb#1{{\fam\msbfam\relax#1}}
\def\widehat#1{\setbox\z@\hbox{$\m@th#1$}%
 \ifdim\wd\z@>\tw@ em\mathaccent"0\msbfam@5B{#1}%
 \else\mathaccent"0362{#1}\fi}
\font\teneufm=eufm10
\font\seveneufm=eufm7
\font\fiveeufm=eufm5
\newfam\eufmfam
\textfont\eufmfam=\teneufm
\scriptfont\eufmfam=\seveneufm
\scriptscriptfont\eufmfam=\fiveeufm
\def\frak#1{{\fam\eufmfam\relax#1}}

\newsymbol\boxdot 1200
\newsymbol\boxplus 1201
\newsymbol\boxtimes 1202
\newsymbol\square 1003
\newsymbol\blacksquare 1004
\newsymbol\centerdot 1205
\newsymbol\lozenge 1006
\newsymbol\blacklozenge 1007
\newsymbol\circlearrowright 1308
\newsymbol\circlearrowleft 1309
\undefine\rightleftharpoons
\newsymbol\rightleftharpoons 130A
\newsymbol\leftrightharpoons 130B
\newsymbol\boxminus 120C
\newsymbol\Vdash 130D
\newsymbol\Vvdash 130E
\newsymbol\vDash 130F
\newsymbol\twoheadrightarrow 1310
\newsymbol\twoheadleftarrow 1311
\newsymbol\leftleftarrows 1312
\newsymbol\rightrightarrows 1313
\newsymbol\upuparrows 1314
\newsymbol\downdownarrows 1315
\newsymbol\upharpoonright 1316
 
\newsymbol\downharpoonright 1317
\newsymbol\upharpoonleft 1318
\newsymbol\downharpoonleft 1319
\newsymbol\rightarrowtail 131A
\newsymbol\leftarrowtail 131B
\newsymbol\leftrightarrows 131C
\newsymbol\rightleftarrows 131D
\newsymbol\Lsh 131E
\newsymbol\Rsh 131F
\newsymbol\rightsquigarrow 1320
\newsymbol\leftrightsquigarrow 1321
\newsymbol\looparrowleft 1322
\newsymbol\looparrowright 1323
\newsymbol\circeq 1324
\newsymbol\succsim 1325
\newsymbol\gtrsim 1326
\newsymbol\gtrapprox 1327
\newsymbol\multimap 1328
\newsymbol\therefore 1329
\newsymbol\because 132A
\newsymbol\doteqdot 132B
 
\newsymbol\triangleq 132C
\newsymbol\precsim 132D
\newsymbol\lesssim 132E
\newsymbol\lessapprox 132F
\newsymbol\eqslantless 1330
\newsymbol\eqslantgtr 1331
\newsymbol\curlyeqprec 1332
\newsymbol\curlyeqsucc 1333
\newsymbol\preccurlyeq 1334
\newsymbol\leqq 1335
\newsymbol\leqslant 1336
\newsymbol\lessgtr 1337
\newsymbol\backprime 1038
\newsymbol\risingdotseq 133A
\newsymbol\fallingdotseq 133B
\newsymbol\succcurlyeq 133C
\newsymbol\geqq 133D
\newsymbol\geqslant 133E
\newsymbol\gtrless 133F
\newsymbol\sqsubset 1340
\newsymbol\sqsupset 1341
\newsymbol\vartriangleright 1342
\newsymbol\vartriangleleft 1343
\newsymbol\trianglerighteq 1344
\newsymbol\trianglelefteq 1345
\newsymbol\bigstar 1046
\newsymbol\between 1347
\newsymbol\blacktriangledown 1048
\newsymbol\blacktriangleright 1349
\newsymbol\blacktriangleleft 134A
\newsymbol\vartriangle 134D
\newsymbol\blacktriangle 104E
\newsymbol\triangledown 104F
\newsymbol\eqcirc 1350
\newsymbol\lesseqgtr 1351
\newsymbol\gtreqless 1352
\newsymbol\lesseqqgtr 1353
\newsymbol\gtreqqless 1354
\newsymbol\Rrightarrow 1356
\newsymbol\Lleftarrow 1357
\newsymbol\veebar 1259
\newsymbol\barwedge 125A
\newsymbol\doublebarwedge 125B
\undefine\angle
\newsymbol\angle 105C
\newsymbol\measuredangle 105D
\newsymbol\sphericalangle 105E
\newsymbol\varpropto 135F
\newsymbol\smallsmile 1360
\newsymbol\smallfrown 1361
\newsymbol\Subset 1362
\newsymbol\Supset 1363
\newsymbol\Cup 1264
 
\newsymbol\Cap 1265
 
\newsymbol\curlywedge 1266
\newsymbol\curlyvee 1267
\newsymbol\leftthreetimes 1268
\newsymbol\rightthreetimes 1269
\newsymbol\subseteqq 136A
\newsymbol\supseteqq 136B
\newsymbol\bumpeq 136C
\newsymbol\Bumpeq 136D
\newsymbol\lll 136E
 
\newsymbol\ggg 136F
 
\newsymbol\circledS 1073
\newsymbol\pitchfork 1374
\newsymbol\dotplus 1275
\newsymbol\backsim 1376
\newsymbol\backsimeq 1377
\newsymbol\complement 107B
\newsymbol\intercal 127C
\newsymbol\circledcirc 127D
\newsymbol\circledast 127E
\newsymbol\circleddash 127F
\newsymbol\lvertneqq 2300
\newsymbol\gvertneqq 2301
\newsymbol\nleq 2302
\newsymbol\ngeq 2303
\newsymbol\nless 2304
\newsymbol\ngtr 2305
\newsymbol\nprec 2306
\newsymbol\nsucc 2307
\newsymbol\lneqq 2308
\newsymbol\gneqq 2309
\newsymbol\nleqslant 230A
\newsymbol\ngeqslant 230B
\newsymbol\lneq 230C
\newsymbol\gneq 230D
\newsymbol\npreceq 230E
\newsymbol\nsucceq 230F
\newsymbol\precnsim 2310
\newsymbol\succnsim 2311
\newsymbol\lnsim 2312
\newsymbol\gnsim 2313
\newsymbol\nleqq 2314
\newsymbol\ngeqq 2315
\newsymbol\precneqq 2316
\newsymbol\succneqq 2317
\newsymbol\precnapprox 2318
\newsymbol\succnapprox 2319
\newsymbol\lnapprox 231A
\newsymbol\gnapprox 231B
\newsymbol\nsim 231C
\newsymbol\ncong 231D
\newsymbol\diagup 231E
\newsymbol\diagdown 231F
\newsymbol\varsubsetneq 2320
\newsymbol\varsupsetneq 2321
\newsymbol\nsubseteqq 2322
\newsymbol\nsupseteqq 2323
\newsymbol\subsetneqq 2324
\newsymbol\supsetneqq 2325
\newsymbol\varsubsetneqq 2326
\newsymbol\varsupsetneqq 2327
\newsymbol\subsetneq 2328
\newsymbol\supsetneq 2329
\newsymbol\nsubseteq 232A
\newsymbol\nsupseteq 232B
\newsymbol\nparallel 232C
\newsymbol\nmid 232D
\newsymbol\nshortmid 232E
\newsymbol\nshortparallel 232F
\newsymbol\nvdash 2330
\newsymbol\nVdash 2331
\newsymbol\nvDash 2332
\newsymbol\nVDash 2333
\newsymbol\ntrianglerighteq 2334
\newsymbol\ntrianglelefteq 2335
\newsymbol\ntriangleleft 2336
\newsymbol\ntriangleright 2337
\newsymbol\nleftarrow 2338
\newsymbol\nrightarrow 2339
\newsymbol\nLeftarrow 233A
\newsymbol\nRightarrow 233B
\newsymbol\nLeftrightarrow 233C
\newsymbol\nleftrightarrow 233D
\newsymbol\divideontimes 223E
\newsymbol\varnothing 203F
\newsymbol\nexists 2040
\newsymbol\Finv 2060
\newsymbol\Game 2061
\newsymbol\mho 2066
\newsymbol\eth 2067
\newsymbol\eqsim 2368
\newsymbol\beth 2069
\newsymbol\gimel 206A
\newsymbol\daleth 206B
\newsymbol\lessdot 236C
\newsymbol\gtrdot 236D
\newsymbol\ltimes 226E
\newsymbol\rtimes 226F
\newsymbol\shortmid 2370
\newsymbol\shortparallel 2371
\newsymbol\smallsetminus 2272
\newsymbol\thicksim 2373
\newsymbol\thickapprox 2374
\newsymbol\approxeq 2375
\newsymbol\succapprox 2376
\newsymbol\precapprox 2377
\newsymbol\curvearrowleft 2378
\newsymbol\curvearrowright 2379
\newsymbol\digamma 207A
\newsymbol\varkappa 207B
\newsymbol\Bbbk 207C
\newsymbol\hslash 207D
\undefine\hbar
\newsymbol\hbar 207E
\newsymbol\backepsilon 237F

\catcode`\@=12

\magnification=\magstep1
\font\title = cmr10 scaled \magstep2

\newcount\refcount
\newcount\seccount
\newcount\sscount
\newcount\eqcount
\newcount\boxcount
\newcount\testcount
\newcount\bibcount
\boxcount = 128
\seccount = -1
\def\sec#1{\advance\seccount by 1\bigskip\goodbreak\noindent
	{\bf\S\number\seccount.\ #1}\medskip \sscount = 0\eqcount = 0}
\def\ss{\advance\sscount by 1\medskip\goodbreak\noindent
	{\bf{(\number\seccount.\number\sscount)}\ \ } \eqcount = 0}
\def\remark#1{\advance\sscount by 1\medskip\goodbreak\noindent{\it #1}
	{\bf{(\number\seccount.\number\sscount)}}\ \eqcount = 0}
\def\proc#1#2{\advance\sscount by 1\eqcount = 0
	\medskip\goodbreak\noindent{\it #1}
	{\bf{(\number\seccount.\number\sscount)}}\ \ {\it #2}}
\def\nproc#1#2#3{\advance\sscount by 1\eqcount = 0\global
	\edef#1{(\number\seccount.\number\sscount)}	
	\medskip\goodbreak\noindent{\it #2}
	{\bf{(\number\seccount.\number\sscount)}}\ \ {\it #3}}
\def\proof{\medskip\noindent{\bf Proof.\ \ }}
\def\eql#1{\global\advance\eqcount by 1\global
	\edef#1{(\number\seccount.\number\sscount.\number\eqcount)}\leqno{#1}}
\def\ref#1#2{\advance\refcount by 1\global
	\edef#1{[\number\refcount]}\setbox\boxcount=
	\vbox{\item{[\number\refcount]}#2}\advance\boxcount by 1}
\def\biblio{{\frenchspacing
	\bigskip\goodbreak\centerline{\bf REFERENCES}\medskip
	\bibcount = 128\loop\ifnum\testcount < \refcount
	\goodbreak\advance\testcount by 1\box\bibcount
	\advance\bibcount by 1\repeat\medskip}}
\def\tightmatrix#1{\null\,\vcenter{\normalbaselines\mathsurround=0pt
	\ialign{\hfil$##$\hfil&&\ \hfil$##$\hfil\crcr
	\mathstrut\crcr\noalign{\kern-\baselineskip}
	#1\crcr\mathstrut\crcr\noalign{\kern-\baselineskip}}}\,}


\def\char{\chi}
\def\geom{{\rm geom}}
\def\tr{{\rm tr}}
\def\scirc{{\scriptstyle\circ}}
\def\C{{\Bbb C}}
\def\Z{{\Bbb Z}}
\def\Q{{\Bbb Q}}
\def\F{{\Bbb F}}
\def\R{{\Bbb R}}
\def\N{{\Bbb N}}
\def\P{{\Bbb P}}
\def\a{{\bf a}}
\def\Sym{{\rm Sym}}
\def\LP{\left(}
\def\RP{\right)}
\def\LB{\left[}
\def\RB{\right]}
\def\LA{\left<}
\def\RA{\right>}
\def\Fr{F}
\def\t{\otimes}
\def\rt#1{{\buildrel #1\over\longrightarrow}}
\def\cF{{\cal F}}
\def\A{{\Bbb A}}
\def\M{{\cal M}}
\def\g{{\frak g}}
\def\h{{\frak h}}
\def\Spec{{\rm Spec}}
\def\Aut{{\rm Aut}}
\def\Hom{{\rm Hom}}
\def\Lie{{\rm Lie}}
\def\rank{{\rm rk}}
\def\GL{{\rm GL}}
\def\Sp{{\rm Sp}}
\def\SU{{\rm SU}}
\def\SO{{\rm SO}}
\def\pialg{\pi_1^{\rm alg}}

\ref\Bour{N. Bourbaki, Groupes et alg\`ebres de Lie, Paris: Masson (1981)}
\ref\Atlas{J. H. Conway, R. T. Curtis, S. P. Norton, R. A.
	Parker, R. A. Wilson, {\it Atlas of Finite Groups},
	Clarendon Press, Oxford, 1985.}
\ref\CR{C. W. Curtis, I. Reiner, {\it Representation Theory of Finite
	Groups and Associative Algebras}, Wiley-Interscience, London, 1962.}
\ref\Davidoff{G. Davidoff, Distribution properties of certain exponential
	sums, {\it Duke J. of Math.} {\bf 62} No. 1 (1991) 205--235.}
\ref\WeilII{P. Deligne, La Conjecture de Weil, II, {\it Publ. Math.
	IHES}, {\bf 52} (1980), 138--252.}
\ref\DM{P. Deligne, D. Mumford, The irreducibility of the space of curves
	of given genus, {\it Publ. Math. IHES}, {\bf 36} (1969) 75--110.}
\ref\DMOS{P. Deligne, J. S. Milne, A. Ogus, K. Shih, {\it Hodge Cycles,
	Motives, and Shimura Varieties}, Lecture Notes in Math. 900,
	Springer-Verlag, Berline, 1982.}
\ref\Katz{N. Katz, {\it Gauss Sums, Kloosterman Sums, and Monodromy Groups},
	Annals of Math. Studies 116, Princeton University Press,
	Princeton, 1988.}
\ref\KatzII{N. Katz, {\it Exponential Sums and Differential Equations},
	Annals of Math. Studies 124, Princeton University Press,
	Princeton, 1990.}	
\ref\Kumar{S. Kumar, Proof of the Parthasarathy-Ranga Rao-Varadarajan
	conjecture, {\it Invent. Math.} {\bf 93} (1988) 117--130.}
\ref\Milne{J. S. Milne, {\it \'Etale Cohomology}, Princeton University
	Press, Princeton, 1980.}
\ref\Mumford{D. Mumford, {\it Curves and their Jacobians}, University of
	Michigan Press, Ann Arbor, 1976.}
\ref\Weyl{H. Weyl, {\it The Classical Groups}, Princeton University Press,
	Princeton, 1939.}
\centerline{The Normal Distribution as a Limit of Generalized
Sato-Tate Measures}
\bigskip
\centerline{Michael Larsen\footnote*{Supported by N.S.A. Grant No. MDA 904-92-H-3026}}
\centerline{University of Pennsylvania}
\centerline{Philadelphia, PA 19104}
\centerline{\tt larsen@math.upenn.edu}
\bigskip
\sec{Introduction}
Let $X$ denote a curve of genus $g$ over a finite field $\F_q$.
By the Riemann hypothesis for function fields,
$$\vert X(\F_q)\vert = 1 - T + q,$$
where
$$- 2g\sqrt{q}\le T\le 2g\sqrt{q}.\eql\RiemannBound$$
We consider this inequality from a naive probabilistic point of view.
Suppose $g = 1$, for example.  Then $E$ has a model as a projective
cubic curve.
Let $F(x,y,z)$ denote the corresponding homogeneous cubic polynomial.
There are $q^2 + q + 1$ $\F_q$-rays through the origin in
affine $3$-space, and loosely speaking, each ray has a probability
$1/q$ of lying on the zero-locus of $F$.  Therefore, $X(\F_q)$ should
should have about $q+1$ points, with an expected error of $O(\sqrt{q})$.
Thus \RiemannBound\ gives the right order of magnitude.  On the other hand,
rather than satisfying a Gaussian distribution law as one might expect,
$T/\sqrt{q}$ is absolutely bounded by $2$.  More precisely, in the
limit $q\to\infty$,
the values of $T/\sqrt{q}$ are uniformly distributed with respect to
a certain measure $\mu_1$, known as the {\it Sato-Tate measure},
if $E$ is drawn at random from the set of
isomorphism classes of elliptic curves over $\F_q$.
For each value of $g$, we obtain in this way a measure $\mu_g$
supported on
$[-2g,2g]$.  The main result of this paper is that the limit of these
measures is, in fact, $\mu = (2\pi)^{-1/2}e^{-x^2/2}\,dx$.
Heuristically speaking,
for a random curve of random genus, over a random finite field,
$T/\sqrt{q}$ is normally distributed.

The distribution of $|X(\F_q)|$, as $X$ varies over a family of
varieties, can be
approached via the cohomological theory of exponential sums,
due to Deligne \WeilII.
More generally, to a family of exponential sums (suitably defined)
one can associate a compact Lie group $G$, the {\it geometric monodromy group},
and a finite-dimensional representation $(\rho,V)$ of $G$ with
character $\chi$.  Under suitable hypotheses, the distribution of values of
the sums is the same as the distibution of values of $\chi(g)$, as
$g$ is drawn randomly from the uniform (Haar) measure on $G$.
For each genus $g\ge 2$, we construct such a group $G_g$ and
a symplectic representation
$V_g$.  It turns out that $G_g$ is the full (compact) symplectic group
$\Sp(2g)$ and $V_g$ is its standard $2g$-dimensional representation.
It is not easy to explicitly compute the distribution of
$\chi(g)$ as $g$ ranges over $\Sp(2g)$, but the moments of the
distribution can be read off from the invariant theory of $\Sp(2g)$.
There is a precise sense in which the invariant theory ``stabilizes''
for large values of $g$.  The measure-theoretic counterpart of this
fact is the statement that $\mu_g$ converges to
$\mu$.

Of course, the heuristic principle suggesting this theorem is not limited
to curves.  One would also expect, for example, that a ``random'' hypersurface
in $\P^3$ of degree $k\gg 0$
should have $q^2 + q + 1 + T$ points, where $T/q$ is normally
distributed.  The corresponding $V_k$ in this case is orthogonal.
Fortunately, the invariant theory of $\SO(n)$ also stabilizes (though one
should be careful to distinguish even and odd values of $n$), and the
limit of the Sato-Tate measures associated to the standard representations
is again $\mu$.  To prove that the Sato-Tate
measures of $(G_k,V_k)$ for hypersurfaces of degree $k$
do indeed tend $\mu$, one would need to show somehow that $G_k$ is ``as large
as possible.''  In most known instances, geometric monodromy groups
associated with self-dual representations do turn out to be either $\Sp$, $O$,
or $\SO$.  See, for instance, \Katz\ 11.1 and \KatzII\ 8.13
for a discussion of Kloosterman sums and hypergeometric sums respectively.

This paper is organized as follows.  The first section analyzes
the geometric monodromy group $(G'_{2g+2},V'_{2g+2})$
of a subfamily of the curves of genus
$g$, namely, the hyperelliptic curves of degree $2g+2$.
We compute the first $2g+1$
moments of the associated Sato-Tate measure $\mu'_{2g+2}$
explicitly.  These turn out
to be the same as the first $2g+1$ moments of $\mu$.
To conclude that $\mu$ is actually the weak limit of the
$\mu'_{2g+2}$ requires some work.
The second section gives an effective version of the Weierstrass approximation
theory which justifies this conclusion.
At this stage, we still do not know that $G'_{2g+2}$ is $\Sp_{2g}$, though
the values of the moments are highly suggestive.  The third section
is pure invariant theory, devoted to an analysis of the fourth moment
of the Sato-Tate measure of a self-dual representation of an infinite
compact group.
Its value is $3$ when $(G,V)$ consists of the standard representation
of a symplectic, orthogonal, or special orthogonal group, and with one
exception, the converse is also true.  With this result in hand, we
return to $G'_{2g+2}$ in the fourth section and prove that is exactly
$\Sp_{2g}$.  Since the full moduli space of curves, $\M_g$, is larger
than the family of hyperelliptic curves, we expect the
geometric monodromy group $G_g$ for $\M_g$ to
be at least at large as $G'_{2g+2}$.  We make this precise, and prove that
in fact $G_g = G'_{2g+2}$.

The results of this paper are superficially similar to those of
Davidoff \Davidoff.  However, the apparent similarity is illusory.
Let $G_n = U(1)^n$ with its standard $n$-dimensional representation.
If the associated Sato-Tate measures are {\it rescaled} by a factor of
$\sqrt{n}$, by the central limit theorem, they will converge to the normal
distribution.  This is a completely different phenomenon from the
(rather miraculous) fact that the Sato-Tate measures of the standard representations of orthogonal and symplectic groups, without any
renormalization, tend to $\mu$. Davidoff considers a special family of
hyperelliptic curves which are endowed with many endomorphisms, and as
a result, she obtains a commutative monodromy group.

This paper has benefited from discussions with C. Epstein, N. Katz,
A. Kouvidakis, G. Kuperberg, A. Lindenstrauss, and R. Pink.
It gives me pleasure to acknowledge their assistance.

\sec{Hyperelliptic Exponential Sums}

In this section we compute the $k^{\rm th}$ moment of the generalized
Sato-Tate measure of $H^1$ of the universal hyperelliptic curve of genus
$g$ for $k < 2g + 2$.  These measures behave as one would expect if
the geometric monodromy group of the sheaf of cohomology groups were
$\Sp(2g)$.  This cohomological point of view is explained below.
We begin with some estimates of exponential sums.

\ss Let $p$ denote a fixed odd prime and $n\ge 3$ an integer.
Let $R_n/\F_p$ denote
the affine variety
$$R_n = \Spec\LP\F_p\LB t_1,\ldots,t_n,{1\over t_1-t_2},
{1\over t_1-t_3},\ldots, {1\over t_{n-1}-t_n}\RB\RP$$
of ordered $n$-tuples with no repeated elements.
We consider the universal hyperelliptic curve $X$ over
$R_n$ given by the homogeneous equation
$$y^2 z^{n-2} - \prod_{i=1}^n (x - t_i z) = 0.\eql\hyperCurve$$
The unique point at $\infty$ on this curve has projective coordinates
$(0,1,0)$.  Let  $\F_q$ denote a finite extension of
$\F_p$ and $\char$ the extension by zero of the quadratic
character on $\F_q^\times$
For each $n$-tuple $\a = (a_1,\ldots,a_n)\in\F_q^n$, and each
$x\in\F_q$, the cardinality of the set of solutions $y\in\F_q$ of
$$y^2 = (x-a_1)\cdots(x-a_n)$$
is $1 + \char((x-a_1)\cdots(x-a_n))$.  Summing over $x$, we
see that the cardinality
of $X_{\a}(\F_q)$ is $q - T + 1$, where
$$T_{\a} = -\sum_{x\in\F_q}\char((x-a_1)\cdots(x-a_n)).\eql\directImage$$

\nproc\RH
{Lemma}{If $P(x)$ a polynomial over a finite field $\F_q$,  then
$$\Bigm\vert\sum_{x\in\F_q}\char(P(x))\Bigm\vert
\le (\deg(P)-1)\sqrt{q}$$
unless $P(x) = c Q(x)^2$ for some $Q(x)\in\F_q[x]$, $c\in\F_q$.}

\proof If $P(x) = R(x)Q(x)^2$, then
$$\Bigm\vert\sum_{x\in\F_q}\char(P(x)) - \sum_{x\in\F_q}\char(R(x))\Bigm\vert
\le\deg(Q).$$
Therefore, without loss of generality, we may assume that $P(x)$ is squarefree.
Let $d$ denote the degree of $P$.  If $d = 0$, there is nothing to prove.
If $d = 1$ or $d = 2$, the sum is always $0$ or $\pm1$ respectively.
Otherwise $y^2 z^{d-2} = z^d P(x/z)$ defines a projective hyperelliptic
curve $X$ of genus $g = \LB{d-1\over 2}\RB > 0$.  By the Lefschetz
trace formula,
$$|X(\F_q)| = q - \tr(\Fr_q: H^1(X,\Q_\ell)) + 1,$$
where $\Fr_q$ denotes the Frobenius element $x\mapsto x^q$, and
$H^1(X,\Q_\ell)$ the \'etale cohomology of $X$ with $\ell$-adic
coefficients, $\ell\neq p$.  By \directImage,
$-\sum_{x\in\F_q}\char(P(x))$ is the trace of Frobenius on $H^1$.
By the Riemann hypothesis for curves over finite fields, this has absolute
value less than or equal to $2g\sqrt{q}$.

\proc{Lemma}{Let $P(x)$ be a monic polynomial over a finite field $\F_q$,
not a perfect square.  Then
$$\sum_{(a_1,\ldots,a_n)\in R_n(\F_q)} \char(P(a_1)P(a_2)\cdots P(a_n))
= O(q^{n/2}).$$}
\proof For each partition $\pi$ of a set $\Sigma$, we write $V_{\pi,\Sigma}$
for the subvariety of the $|\Sigma|$-dimensional affine space $\A^\Sigma$
defined by the equations $x_a = x_b$ whenever
$\{a,b\}\subset\kappa\in\pi$.  Setting $\Sigma_n = \{1,2,\ldots,n\}$,
we consider
$$S_{\pi,\Sigma_n} =
\sum_{(a_1,\ldots,a_n)\in V_{\pi,\Sigma_n}(\F_q)}
\char(P(a_1)P(a_2)\cdots P(a_n))$$
and prove by induction on $n$ that all of
such sums are $O(q^{n/2})$.  If $\kappa\in\pi$ is a singleton $\{k\}$, then
$$S_{\pi,\Sigma_n} = \LP\sum_{a_k\in\F_q} \char(P(a_k))\RP
S_{\pi\setminus\{\kappa\},\Sigma_n\setminus\{k\}} = O(q^{n/2})$$
by Lemma \RH\ and induction on $n$.  If all the equivalence classes
of $\pi$ have two elements or more, there can be no more than $n/2$ of them,
so $\dim(V_{\pi,\Sigma_n})\le n/2$.  It follows that
$|S_{\pi,\Sigma_n}|\le q^{n/2}$.

By the inclusion-exclusion
principle, the characteristic function on $R_n$ is a linear combination of
the characteristic functions of $V_{\pi,\Sigma_n}$.  The lemma follows.

\nproc\MainEstimate{Proposition}{Let
$$F(m) = \cases{(m-1)!!:=(m-1)(m-3)(m-5)\cdots(5)(3)(1)&if $m$ is even,\cr
0&if $m$ is odd.\cr}$$}
For $m<n$,
$$\sum_{\a\in R_n(\F_q)} (-T_{\a})^m = F(m) q^{n + m/2} + O(q^{n + m/2 - 1}).$$
\proof By \directImage,
$$\sum_{\a\in R_n(\F_q)} (-T_{\a})^m = 
\sum_{\a\in R_n(\F_q)} \sum_{x_1\in\F_q}\cdots\sum_{x_m\in\F_q}
\char\LP\prod_{i=1}^n\prod_{j=1}^m (x_j - a_i)\RP.$$
Defining
$$P_{x_1,\ldots,x_m}(y) = \prod_{j=1}^m (y - x_j) = (-1)^m\prod_{j=1}^m 
(x_j-y),$$
we obtain
$$\sum_{\a\in R_n(\F_q)} (-T_{\a})^m = \sum_{x_1\in\F_q}\cdots\sum_{x_m\in\F_q}
\sum_{\a\in R_n(\F_q)}
\char((-1)^{mn})
\char\LP P_{x_1,\ldots,x_m}(a_1)\cdots P_{x_1,\ldots,x_m}(a_n)\RP.
\eql\switched$$
If there exists a partition $\pi$ of $\Sigma_m$ into pairs such that
$(x_1,\ldots,x_m)\in V_{\pi,\Sigma_m}(\F_q)$, then as a polynomial in $y$,
$P_{x_1,\ldots,x_m}(y)$ is a perfect square,
but otherwise it is not.  Therefore,
outside the union of such loci, the inner sum of \switched\ is
of order $O(q^{n/2})$.  There are $F(m)$ partitions $\pi$ of $\Sigma_m$
into pairs, and the loci $V_\pi$ intersect
pairwise in sets of lower dimension, so
$$\sum_{\a\in R_n(\F_q)}\hskip -10pt T_{\a}^m\!
= (-1)^{mn-m} F(m) q^{n + m/2}\!+ O(q^{n + m/2 - 1}\!+ q^{m + n/2})
= \!F(m) q^{n + m/2}\!+ O(q^{n + m/2 - 1}).$$

\remark{Remark}If $m$ and $n$ are both odd, it is easy to see that
$$\sum_{\a\in R_n(\F_q)} T_{\a}^m = 0.$$
Indeed, if $\lambda\in\F_q$ is chosen such that $\char(\lambda) = -1$,
then
$$T_{\lambda\a} = -\sum_{x\in\F_q}\char((x-\lambda a_1)\cdots(x-\lambda a_n))
= -\sum_{\lambda x\in\F_q}\char(\lambda^n)\char((x- a_1)\cdots(x- a_n))
= - T_{\lambda\a}.$$

\ss Let $S/\F_p$ denote a geometrically connected variety and $\bar s$
a fixed geometric point of $S$.
The functor assigning to each finite \'etale $X/S$
the set $\Hom_S(\bar s,X)$, is represented by a projective system
$\{X_i,\ \phi_{ij}:X_j\to X_i\}$ of finite \'etale covers of $S$.
As usual, we define
$$\pialg(S,\bar s):=\varprojlim_i \Aut_X(X_i).$$
Given a second geometric point $\bar s'$, there is an isomorphism
$\pialg(S,\bar s')\cong\pialg(S,\bar s)$ canonical up to inner automorphism
({\it i.e.}, up to choice of ``path'') \Milne\ I 5.1 (a).
Let $\bar S$ denote the variety obtained from
$S$ by extension of scalars to $\overline{\F}_p$.  There is a short exact
sequence
$$0\to\pialg(\bar S,\bar s)\to\pialg(S,\bar s)\to\widehat{\Z}\to 0$$
(\WeilII\ 1.1.13).
By a smooth $\ell$-adic sheaf $\cF$ on $S$, we mean a continuous
homomorphism $\pialg(S,\bar s)\to\GL_n(\Q_\ell)$.  We fix
an embedding $\iota:\Q_\ell\to\C$ and say
that $\cF$ is {\it pure of weight $w\in\Z$} if for every
$\F_q$-point $x\in S(\F_q)$, every geometric point $\bar x$ of $x$,
and every choice of path $\pialg(S,\bar x)\to\pialg(S,\bar s)$,
the image of Frobenius under the composition
of arrows
$$\pialg(x,\bar x)\to\pialg(S,\bar x)\to\pialg(S,\bar s)\to\GL_n(\Q_\ell)
\to\GL_n(\C)$$
has all its eigenvalues of absolute value $q^{w/2}$.
We write $\tr(\cF_x)$ for the trace of the resulting matrix.

\ss Let $d = \dim(S)$.  By Poincar\'e duality (\Milne\ VI 11.1)
$$H^{2d}_c(\bar S,\cF)\tilde\to \Hom_{\Q_\ell}(H^0(\bar S,\cF^\vee(d)),\Q_\ell)
= \cF^{\pialg(\bar S,\bar s)}(-d).$$
The right hand side is pure of weight $w + 2d$.
Let $G^{\geom}$ denote the Zariski closure of $\pialg(\bar S,\bar s)$
in $\GL_n(\Q_\ell)$.  The space of $G^{\geom}$-invariants in $\Q_\ell^n$
is the same as the space of $\pialg(\bar S,\bar s)$-invariants, so
$$\dim(H^{2d}_c(\bar S,\cF)) = \dim\bigl(\cF^{G^\geom}\bigr).$$
Note that by \WeilII\ 1.3.9, 3.4.1 (iii), the identity
component of $G^{\geom}$ is semisimple.
By the Lefschetz trace formula \Milne\ VI 13.4,
$$\sum_{x\in S(\F_q)} \tr(\cF_x)=
\sum_{i=0}^{2d} (-1)^i\tr(\Fr_q : H^i_c(\bar S,\cF)).$$
By \WeilII\ 3.3.1, the $i = 2d$ term dominates the sum for $q\gg 0$.
Therefore, if
$$a_n = p^{n(-d-w/2)}\sum_{x\in S(\F_{p^n})} \tr(\cF_x)$$
tends to a limit as $n\to\infty$, then $H^{2d}_c(\bar S,\cF)$ must
be just $\Q_\ell(-d-w/2)$, and
$$\lim_{n\to\infty} a_n = \dim\bigl(\cF^{G^\geom}\bigr).\eql\limitFormula$$

\nproc\InvariantDim
{Proposition}{Let $X\rt{\pi}R_n$ denote the universal hyperelliptic
curve over the base $R_n$, and $\cF = R^1\pi_*(\Q_\ell)$ the sheaf of
$H^1$ of the fibres over $R_n$.
Then $\cF$ is pure of weight $1$ and for $m < n$,
$$\dim\LP\LP\cF^{\t m}\RP^{G^\geom}\RP = F(m),$$
where $F(m)$ is defined as in Prop. \MainEstimate.}

\proof The purity statement is immediate from \WeilII\ 3.4.11.
This implies that $\cF^{\t m}$ is pure of weight $m$ for all positive
integers $m$.  The proposition now follows immediately from Prop.
\MainEstimate\ and \limitFormula.

\ss Let $G$ denote an algebraic group over $\C$
with semisimple identity component and $(\rho,V)$ a complex representation
of $G$.  There is a unique compact real form
$G^c$ of $G$, and
$$\dim(V^G) = \dim\bigl(V^{G^c}\bigr).$$
On the other hand, if $\mu$ denotes Haar measure on $G^c$, then
$$\dim\bigl(V^{G^c}\bigr) = \int_{G^c}\chi_\rho(g)\mu,\quad\chi_\rho
:= \tr\scirc\rho.$$
We define the {\it generalized Sato-Tate} measure for $(G^c,V)$
as the direct image measure
$$\mu_{S-T}:={\chi_\rho}_*\mu$$
on $\R$.  In particular, we can associate a Sato-Tate measure to
every pure sheaf over $S$.
By definition of direct image,
$$\int_{-\infty}^\infty x^n\mu_{S-T} = \int_{G^c}\chi_\rho(g)^n\mu
= \int_{G^c}\chi_{\rho^{\t n}}(g)\mu = \dim\bigl({V^{\t n}}^G\bigr).
\eql\momentDim$$
Prop. \InvariantDim\ can be interpreted as the computation of the
first $n-1$ moments of the Sato-Tate measure associated with the family
of hyperelliptic curves over $R_n$.

\sec{Moments and Limits of Non-negative Measures}

\ss Let $f:\R\to\R^{> 0}$ denote a positive real valued-function.
We say $f$ is {\it very rapidly decreasing} (VRD) if
every smooth compactly supported function $g$ is a uniform limit
(on~$\R$) of functions of the form $f(x)P(x)$, where $P(x)$ is a
polynomial.

\proc{Lemma}{If $f(x)$ is VRD, then so is $f(\lambda x)$ for all $x$.}
\proof The ring of polynomials and the $L^\infty$ norm are both invariant
under the rescaling $x\mapsto\lambda x$.

\proc{Lemma}{If $f(x)$ is VRD and $g(x)$ is positive, bounded continuous
function, then $f(x)g(x)$ is VRD.}

\proof Approximating a smooth compactly supported function $h(x)$
by $f(x)g(x)P(x)$ in the $L^\infty$ topology is equivalent to approximating
$h(x)/g(x)$ by $f(x)P(x)$.

\proc{Lemma}{If $f(x)$ is VRD, then the moment integrals
$$\int_{-\infty}^\infty f(x)x^n\,dx$$
converge for all $n$.}

\proof The space of smooth compactly supported functions on $\R$ is of infinite
dimension, and all such spaces are uniform limits of sequences $f(x) P_n(x)$.
Therefore, $f(x) P(x)\in L^\infty(\R)$ for some polynomial $P$ of arbitrarily
large degree.  Therefore $f(x) = o(x^{-n})$ for all $n$, and the lemma follows.

\proc{Proposition}{The normal functions $f(x) = e^{-Kx^2}$ are VRD.}

\proof Without loss of generality we may assume $K = 1/2$.
Let $T_n(x)$ denote the $n^{\rm th}$ Chebyshev polynomial
of the first kind.  When $x\in[-1,1]$ and $4\vert n$, we have
$$T_n(x) = \cos(n\cos^{-1}x) = \cos(n\sin^{-1}x).$$
For $|x|\le 1 < {\pi\over 3}$,
$$1 - {x^2\over 4}\ge \cos x\ge 1 - {x^2\over 2},\ {x^2\over 4}\le
\sin^2 x\le x^2,$$
so when $4\vert n$,
$$1 - 2n^2x^2\le  T_n(x)\le 1 - {n^2\over 4}x^2$$
for $|x|\le {1\over 2n}$.
Let $m$ denote an integer congruent to $2$ (mod $4$).  For $|x| < m/2$,
$$1 - {2x^2\over m^2}\le 1 - 2(m-2)^2\LP{x\over 2m}\RP^2
\le T_{m-2}\LP{x\over m^2}\RP\le 1 - {(m-2)^2\over 4}\LP{x\over m^2}\RP^2
\le 1 - {x^2\over 16 m^2}.\eql\innerBound$$

We define
$$Q_m(x) = \LB\LP 1-{x^2\over m^4}\RP T_{m-2}\LP{x\over m^2}\RP\RB^{m^3/4}.$$
As $(1-\epsilon)^{1/\epsilon}\ge 1/4$ for $0\le \epsilon\le {1\over 2}$,
for $|x| < m^{-1/2}$,
$$Q_m(x) \ge ((1 - m^{-5})(1 - 2m^{-3}))^{m^3/4} \ge (1 - 4m^{-3})^{m^3/4}
\ge {1\over 4}.$$
As $m^3/4$ is even, $Q_m(x)\ge 0$ for all $x$, so
$$I_m = \int_{-m^2}^{m^2} Q_m(x)\ge {1\over 2\sqrt{m}}.$$
In the other direction, \innerBound\ implies
$$Q_m(x)\le\LP1-{x^2\over 16m^2}\RP^{m^3/4}\le e^{-mx^2/64}\eql\secondRing$$
for $|x| < m/2$.  For $m/2 \le |x| \le \sqrt{2m^4 - m^2/4}$,
$$Q_m(x)\le \LP1-{x^2\over m^4}\RP^{m^3/4}\le e^{-m/16}.\eql\thirdRing$$
Finally, we note that the identity $T_n(\cosh x) = \cosh nx$
implies that for $|y|\ge 1$,
$$T_n(y) < |2y|^n.$$
For $|x|\ge m^2$, then,
$$0 < Q_m(x) e^{-(x-r)^2/2} < (2^m m^{-2m} x^m)^{m^3/4} e^{-(x-r)^2/2}.$$
From the fact that $|x^n e^{-x^2/2}|$ achieves its maximum when
$x = \pm\sqrt{n}$, we deduce
$$\sup_{x\in\R} x^n e^{-(x-r)^2/2} = \sup_{x\in\R} (x+r)^n e^{-x^2/2}
= \sup_{x\in\R}(1+r/x)^n e^{-x^2/2}.$$
If the maximum is achieved for some $x \ge 2r$, then
$$\sup_{x\in\R} x^n e^{-(x-r)^2/2} \le n^{n/2} C^{-n/2},$$
where $C = 4e/9 > 1$.  Thus, if $r < \sqrt{n\over 9C}$,
$$\sup_{x\in\R} x^n e^{-(x-r)^2/2} \le \sup(\sup_{|x|\le 2r}x^n e^{-(x-r)^2/2},
n^{n/2} C^{-n/2}) \le n^{n/2} C^{-n/2}.$$
We conclude that if $r < {m^2\over 6\sqrt{C}}$,
$$Q_m(x) e^{-(x-r)^2/2} < (2^m m^{-2m})^{m^3/4} \sup_{x\in\R} x^{m^4/4}e^{-(x-r)^2/2}
\le C^{-{m^4\over 4}}.\eql\outerRing$$

For each $n\in\N$, we set $m = 4n+2$ and
define $f_n(x)$ as the product of
$Q_m/I_m$ and the characteristic function of the interval
$[-m^2,m^2]$.  Then $f_n$ is a non-negative measure of integral $1$,
and the estimates \secondRing\ and \thirdRing\ imply that
$f_n$ is an approximate identity.  Therefore, if $g$ is a smooth compactly
supported function, the sequence of convolutions $f_n\ast g$ converges
uniformly to $g$.  By \thirdRing\ and \outerRing, since $g$ is supported
on a subset of $[-r,r]$ for some $r$,
$$e^{-x^2/2}((f_n - Q_{4n+2}/I_{4n+2})\ast g)$$
tends to zero uniformly on $\R$ as $n\to\infty$.
Therefore, $g$ is the uniform limit of the products $e^{-x^2/2}((Q_{4n+2}/I_{4n+2})\ast g)$ of a fixed Gaussian with a sequence
of polynomials.

\nproc\LimitMeasure
{Proposition}{Let $\mu_1,\mu_2,\ldots$ denote a sequence of non-negative
measures on the real line and $f(x)$ a smooth
positive real-valued VRD function.  Suppose that for all $m < n$
$$\int_{-\infty}^\infty x^m\mu_n = \int_{-\infty}^\infty x^m f(x)^4\,dx.
\eql\equalMoments$$
Then $\mu_i$ converges to $f(x)^4\,dx$ in the weak-$\ast$ topology.}

\proof We note that $f(x)^4$ is VRD, so the right hand side of
\equalMoments\ converges.
Let $g(x)$ by any smooth compactly supported function which takes only
non-negative values.  Choose $\epsilon > 0$, and let $h(x)$
denote a smooth compactly supported function
such that $h(x) = \sqrt{g(x)+\epsilon}$ on the support of $g$
and $|h(x)^2 - g(x)|_\infty = \epsilon$.
Now $f(x)h(x)$ is the uniform limit of a sequence $f(x)A_n(x)$,
where $A_n(x)\in\R[x]$.  For $n$ sufficiently large
$$(f(x)A_n(x))^2 \ge f(x)^2 g(x)$$
on $\R$.  Letting $\epsilon$ tend to $0$, we can write
$f(x)^2 g(x)$ as a uniform limit of functions $f(x)^2 B_n(x)$,
where $B_n(x)\in\R[x]$, and $B_n(x) \ge g(x)$ for all $x\in\R$.
As $f(x)^2$ is VRD, it is integrable, so
$f(x)^4 B_n(x)$ converges in the $L^1$ norm to $f(x)^4 g(x)$.  Therefore
$$\lim_{n\to\infty}\int_{-\infty}^\infty g(x)\mu_n
\le\varliminf_{m\to\infty}\lim_{n\to\infty}\int_{-\infty}^\infty B_m(x)\mu_n
=\varliminf_{m\to\infty}\int_{-\infty}^\infty B_m(x)\mu
=\int_{-\infty}^\infty g(x)\mu.\eql\fundIneq$$
On the other hand, the measure of the real line is the same with
respect to $\mu$ and with respect to $\mu_n$ for $n\gg 0$, so we must
have equality in \fundIneq.  Finally, every smooth compactly
supported function is the difference of two such functions which
are everywhere non-negative, and the proposition follows.

\proc{Corollary}{If $\mu'_n$ denotes the Sato-Tate measure associated with
the family of hyperelliptic curves over $R_n$, then in the weak-$\ast$
topology,
$$\lim_{n\to\infty}\mu'_n = {1\over\sqrt{2\pi}}e^{-x^2/2}\,dx.
\eql\normalLimit$$}

\proof By Prop. \InvariantDim\ and \momentDim,
$$\int_{-\infty}^\infty x^m \mu'_n = F(m),$$
for $m < n$.   Applying Prop. \LimitMeasure\ to $(2\pi)^{-1/8}e^{-x^2/8}$,
we deduce \normalLimit\ from the integral
$$\int_{-\infty}^\infty x^m e^{-x^2/2}\,dx = \sqrt{2\pi}F(m),$$
which is easily checked by integration by parts.

\sec{A Problem in Invariant Theory}

Throughout this section $G$ will always denote a compact Lie group and
$V$ a faithful finite dimensional representation.

\ss Given a fixed $G$ and $V$, we define the sequence of Sato-Tate moments
$$a_n = a_n(G,V) = \dim\LP{V^{\t n}}^G\RP,\quad n\ge 1.$$
It is possible, in general, that $a_1 = a_2 = \cdots = a_m = 0$,
for any desired value of $m$.  This is the case, for instance, if
$G = \SU(n)$, $V$ is the standard $n$-dimensional representation,
and $n > m$.  However, if $V$ is self-dual, then $a_{2k} > 0$
for all $k$.  Indeed,
$$a_{2k} = \dim\LP{V^{\t 2k}}^G\RP = \dim\LP\Hom_G\LP V^{\t k},V^{\t k}\RP\RP
= \sum_i m_i^2,$$
where $m_i$ denote the multiplicities of the irreducible factors of
$V^{\t k}$.  In particular, $a_{2k}$ is at least as large as the
number of summands appearing when $V^{\t k}$ is decomposed into irreducible
representations.

Henceforth, we shall always assume $V$ is self-dual.

%
\ss This section is devoted to the classification of pairs $(G,V)$ such
that $a_4(G,V) = 3$.  The complete classification problem seems to be quite
difficult.  A wide variety of interesting finite groups admit such representations.  For example, the irreducible
$3$-dimensional representations of $A_5$ satisfy this condition.
So do the standard $6$, $7$, and $8$ dimensional representations
of the Weyl groups of $E_6$, $E_7$, and $E_8$ respectively.
The natural $24$-dimensional representation of the automorphism group
of the Leech lattice does as well.
The 133-dimensional representation of
the Harada-Norton group and the 248-dimensional representation of
the Thompson group provide even more exotic examples.
The equality $a_4 = 3$
can be readily checked in all these cases with the aid of character
tables, such as those in \Atlas.

Fortunately, in our application, finite groups may be ruled
out on geometric grounds, so
we are not obliged to attempt to classify the solutions.  It seems likely
that a complete list is attainable, using the classification
of finite simple groups and the following lemma:

\proc{Lemma}{Let $(G,V)$ denote a solution to the equation $a_4(G,V) = 3$.
Then every for every non-abelian normal subgroup $H$ of $G$ the restriction
of $V$ to $H$ is irreducible.  In particular, the centralizer of $H$ in
$G$ has order $\le 2$.}

\proof First we observe that $V$ must be an irreducible $G$-module.
Indeed, if
$V = V'\oplus V''$, then the trivial representation appears with multiplicity
$\ge 2$ in $V^{\t 2}$.  Therefore, $a_4\ge 4$, contrary to hypothesis.
Let $H$ be a normal subgroup of $G$.  As a representation of $H$,
$V$ decomposes into a direct sum 
$$(W_1\oplus\cdots W_k)\t\C^\ell,\quad
\dim(W_1) = \cdots = \dim(W_k) = m,$$
where $H$ acts irreducibly on the
$W_i$ and trivially on $\C^\ell$ (\CR\ 49.7).
As $H$ is non-abelian and $V$ is a faithful $H$-module, $m > 1$.
Every element of $G$ maps an $H$-isotypic factor $W_i\t\C^\ell$
into another such factor, $W_j\t\C^\ell$.  Therefore,
$$V\t V = \bigoplus_{i=1}^k \Sym^2(W_i\t\C^\ell)\ \oplus
\ \bigoplus_{i=1}^k \Lambda^2(W_i\t\C^\ell)\ \oplus
\ \bigoplus_{i\neq j} (W_i\t\C^\ell)\t(W_j\t\C^\ell)$$
represents $V\t V$ as a direct sum of three $G$-modules.  If
$k>1$ or $\ell > 1$, all three pieces are have dimension $>1$.  Since
there is also at least one $G$-invariant in $V\t V$, $a_4\ge 4$.
Therefore, $k=\ell=1$, and $V$ is an irreducible $H$-module.  By
Schur's lemma, the centralizer of $\rho(H)$ in $\GL(V)$ consists of the
scalar matrices.  But all traces of elements in $\rho(G)$ are real,
so only scalar matrices $\pm1$ are possible.  As $V$ is
faithful on $G$, the centralizer of $H$ has order $\le 2$.

\proc{Corollary}{If $a_4(G,V) = 3$, then the identity component of $G$
is either a torus or a semisimple group.}

\proof Let $G^\circ$ denote the identity component of $G$.
It is a normal subgroup of $G$.
If it is not a torus, then it has a finite centralizer in $G$, hence
a finite center.  Therefore, it is semisimple.

\nproc\Classification
{Proposition}{If $a_4(G,V) = 3$ for an infinite compact group $G$,
then $G\subset\GL(V)$ is $N_{\SU(2)}U(1)\subset U(2)$,
$SO(n)\subset U(n)$, $O(n)\subset U(n)$, or $\Sp(2n)\subset U(2n)$.}

\proof Suppose first that $G^\circ$ is a torus $T$.  Then $V$ is a direct
sum $\chi_1\oplus\cdots\chi_n$ of characters of $T$, and for each
$\chi_i$ there exists $\chi_j = \chi_i^{-1}$.  As $V$ is irreducible
as a $G$-representation, all the characters $\chi_i$ must lie in a single
orbit under the action of $G/T$ on the character group $X^*(T)$.
If some $\chi_i$ were trivial, then all $\chi_i$ would be trivial,
contrary to the assumption that $V$ is faithful.  It follows that (suitably
renumbering the indices), $V$ is the direct sum
$$\chi_1\oplus\chi_1^{-1}\oplus\chi_2\oplus\cdots\oplus
\chi_{n/2}^{-1}.$$
Therefore, as $G$-module, $V\t V$ decomposes into the following three pieces:
$$\bigoplus_i
\LP\chi_i\t\chi_i\oplus\chi_i^{-1}\oplus\chi_i^{-1}\RP
\ \oplus\ \bigoplus_i
\LP\chi_i\t\chi_i^{-1}\oplus\chi_i^{-1}\t\chi_i\RP
\ \oplus\ \bigoplus_{i\neq j}
\LP\chi_i\t\chi_j^{-1}\oplus\chi_i^{-1}\t\chi_j\RP.$$
If $n/2 > 1$, then these pieces are all of dimension $\ge 2$.  Since
$V^{\t 2}$ also has a $G$-invariant line, this implies $a_4\ge 4$.
We conclude that $n=2$, so $G$ is contained in $\SU(2)$.  As $T$
is normal in $G$, $G$ can only be the normalizer of a maximal torus.

Suppose, on the contrary, that $G^\circ$ is a semisimple group.
As $G^\circ$ is normal in $G$, the restriction of an irreducible
representation of $G$ to $G^\circ$ is the direct sum of highest
weight modules $V_{\lambda_i}$ of $G^\circ$, where the $\lambda_i$
lie in the same orbit of the automorphism group $\Gamma$ of the root system
$\Phi$ of
$G^\circ$.  We apply this observation to $V_\lambda^{\t 2}$.
The dual of the Killing form gives a $\Gamma$-invariant
inner product on the space of characters of $G$, so it suffices
to find submodules $V_\mu$ of $V_\lambda^{\t 2}$ of four different
lengths.

Given a semisimple Lie algebra $\g$ and a representation $V$, let
$S(\g,V)$ denote the set of norms $\Vert\lambda\Vert^2$, where
$V_\lambda$ is a submodule of $V$.  If $(\h,W)$ is a second pair, then
$$S(\g\times\h,V\boxtimes W)=\{x+y\mid x\in S(\g,V),\,y\in S(\h,W)\},$$
so
$$\vert S(\g\times\h,V\boxtimes W)\vert\ge \vert S(\g,V)\vert
+\vert S(\h,W)\vert - 1.$$
If $\g$ is simple and $V_\lambda$ is faithful and self-dual, then
$\vert S(\g,V_\lambda)\vert\ge 2$ because $0,\,2\lambda\in S(\g,V_\lambda)$.
By \Kumar,
$V_\lambda^{\t 2}$ contains a submodule $V_\mu$ in every Weyl orbit
$W\mu = W(\lambda + w\lambda)$, for fixed $w\in W$.
If $\rank(\g) > 1$, this implies that $\vert S(\g,V_\lambda)\ge 3$
because $W$ acts irreducibly on the root space, and therefore
some $w \lambda\not\in\{\pm\lambda\}$.
Since $V$ is self dual, $W\lambda$ is
invariant under multiplication by $-1$.
Suppose there exists $w\in W$ such that
$$w\lambda\not\in\{\lambda,-\lambda\}\cup \lambda^\perp.$$
Then $\lambda+\lambda$, $\lambda + w\lambda$, $\lambda - w\lambda$,
and $\lambda - \lambda$ are all of different lengths, so
$\vert S(\g,V_\lambda)\vert\ge 4$.

We have seen that if $\g = \Lie(G^\circ)$, and $V$ is a representation of
$G$ such that $a_4(G,V) = 3$, then $\vert S(\g,V)\vert\le 3$.
In view of the foregoing analysis, this implies that
$\g = {\frak sl}_2 \times {\frak sl}_2$ with $V$ the exterior tensor
product of the two standard representations, or $\g$ is simple
with $V = V_\lambda$, where
$\mu_1,\,\mu_2\in W\lambda$ implies $\mu_1=\pm\mu_2$ or $\mu_1\perp\mu_2$.
By the classification of simple Lie algebras, the latter
condition implies that
$\g$ is of type $A_1$, $B_n$ ($n\ge 2$),
$C_n$ ($n\ge 2$), or $D_n$ ($n\ge 3$), and $\lambda$ is a positive integral
multiple of the fundamental weight $\omega_1 = (1,0,\ldots,0)$,
in the notation of
\Bour\ VI {\it Planches}.
Note that the Lie algebras $B_2$ and $C_2$ are the same,
but the value of $\omega_1$ depends on which name we choose.

Assume that $n = \rank(\g) > 1$.  We recall the Freudenthal formula for the multiplicity $m_\omega(\lambda)$
of a weight $\lambda$ appearing in the highest weight module
$V_\omega$ (\Bour\ VIII \S9 Ex. 5 (g)):
$$m_\omega(\lambda) = {2\sum_{\alpha\in\Phi^+}\sum_{i=1}^\infty
m_\omega(\lambda + i\alpha)\LA\mu+i\alpha,\alpha\RA\over
\LA\omega+\rho,\omega+\rho\RA-\LA\lambda+\rho,\lambda+\rho\RA},$$
where $\rho$ denotes the half sum of roots.
Applying this formula for $B_n$, $C_n$, and $D_n$ to
$\omega = k\omega_1$ and
$\lambda = (k-1,1,0,\ldots,0)$, in each case we get a multiplicity of $1$.
As $(2k-1,1,0,\ldots,0)$ appears with multiplicity $2$ in
$V_{k\omega_1}^{\t 2}$, $V_{(2k-1,1,\ldots,0)}$ is a submodule of
$V_{k\omega_1}^{\t 2}$.  By \Kumar, we know that
there are also submodules $V_{(2k,0,\ldots,0)}$, $V_{(k,k,0,\ldots,0)}$,
and $V_0$.  If $k > 1$, the highest weights of these modules
are all of different lengths.
We conclude that $\lambda = \omega_1$.  As $V$ is faithful,
this means that $G^\circ$ is $SO(2n+1)$, the compact form of $Sp(2n)$,
or $SO(2n)$, depending on whether $\g$ is of type $B_n$, $C_n$, or $D_n$.
In any case, $G$ is contained in the normalizer of the image of
$G^\circ$ under its standard representation, so $G = G^\circ$
in the symplectic case and $G$ can be either of type $SO$ or type $O$
in the orthogonal case.

Finally we consider the ${\frak sl}_2$ cases.  When $\g = {\frak sl}_2$,
$V = V_{k\omega_1}$, then $\vert S(\g,V^{\t 2})\vert = k + 1$.
When $k = 1$, we have the standard representation of $G^\circ = \SU(2)$,
so $a_4(G,V) < 3$.  When $k = 2$, we have the solutions $\SO(3)$
and $O(3)$ enumerated above.  Finally, if
$\g = {\frak sl}_2\times {\frak sl}_2$,
the condition $\vert S(\g,V)\vert \le 3$ implies that $V$ is the
exterior tensor product of the standard representations of the two
factors.  This gives rise to the solutions $\SO(4)$ and $O(4)$ enumerated
above. 

\sec{Monodromy for the Moduli Spaces of Curves}

\proc{Theorem}{For $n\ge 5$, the geometric monodromy $(G'_n, V_n)$ of
$R^1\pi_*\Q_\ell$ for the family of hyperelliptic curves
over $R_n$ is $Sp_{2[(n-1)/2]}$.}

\proof Applying Prop. \MainEstimate\ for $m = 4$, we see that
$a_4(G'_n,V_n) = 3$.  Now $\dim(V_n) = 2g = 2[(n-1)/2]$, so by
Prop. \Classification, $G'_n$ is finite, symplectic, or orthogonal.
On the other hand, the cup product on $H^1$ is anti-symmetric, so by
Poincar\'e duality, $V_n$ is a symplectic representation.
Therefore, it suffices to prove that $G'_n$ is infinite.  Now
$$\bar R_n = \Spec\LP\overline{\F}\LB t_1,\ldots,t_n,{1\over t_1-t_2},
{1\over t_1-t_3},\ldots, {1\over t_{n-1}-t_n}\RB\RP$$
is an open subvariety of affine space $\A^n$ over $\overline{\F}$,
and $X$ extends naturally to the projective curve over $\pi:X\to\A^n$
defined by \hyperCurve.
For every non-singular
curve $Z\subset \A^n$, we define $Z^\circ = Z\cap R_n$.
We choose $Z$ such that $Z\setminus Z^\circ$ is a non-empty subset of
the smooth locus of $\A^n\setminus R_n$.
The restriction $X_Z := X\times_{\A^n} Z$
is a Lefschetz pencil, and the restriction of $R^1\pi_{Z*}\Q_\ell$ to
$Z^\circ$ is a smooth $\ell$-adic sheaf $\cF$ with finite monodromy.
Therefore, there exists a finite \'etale cover $Y^\circ$ of $Z^\circ$
on which $\cF$ is a sheaf with trivial monodromy.  The normalization of
$Y^\circ$ over $Z$ is a non-singular curve $Y$ with a Lefschetz pencil
$X_Y\to Y$.  The fibres over the (non-empty) set $Y\setminus Y^\circ$ have
a double point, so there is at least one vanishing cycle of
$R^1\pi_{Y*}\Q_\ell$.  By the Picard-Lefschetz theorem (\Milne\ V 3.15),
the space of $\pi_1$-invariants of $R^1\pi_{Y^\circ*}\Q_\ell$ is orthogonal
to the space of vanishing cycles under the Poincar\'e pairing.  As the
Poincar\'e pairing is perfect, we have a contradiction, and the theorem holds.

\remark{Remark} In the light of this result, the moment computations
of Prop. \InvariantDim\ give
$$a_m(\Sp(2g),{\rm Std}) =
\cases{0&if $m$ is odd,\cr (m-1)!!&if $m$ is even,\cr}$$
for $m\le 2g+1$.  This is a classical result of invariant theory
\Weyl\ 6.1.A, 6.1.B.  The analogous result also holds for orthogonal
groups \Weyl\ 2.11.A, 2.17.A.

\ss For each pair $(g,n)$ of non-negative integers, consider the functor of $n$-pointed curves of genus $g$ over $\F_p$,
{\it i.e.}, the functor assigning
to each scheme $S/\F_p$ the set of $n+1$-tuples $(\pi,\sigma_1,\ldots,\sigma_n)$, 
where $\pi:X\to S$ is a proper smooth morphism and $\sigma_i:S\to X$ are
sections of $\pi$ such that for all geometric points $\bar s$ on $S$,
$X_{\bar s}$ is a curve of genus $g$ and $\sigma_i(\bar s)$ are
distinct points.  For fixed $g$, when $n$ is sufficiently large,
this functors is represented by a quasi-projective variety $\M_{g,n}$
\Mumford\ II.  By a theorem of Deligne and Mumford \DM,
the coarse moduli space of
curves $M_g$ and therefore the variety 
$\M_{g,n}$ is geometrically irreducible.
Let $Y_{2g+2,n}$ denote the complement of the diagonal on the
$n$-fold fibre power of the universal curve on $R_{2g+2}$.  From the projection
$Y_{2g+2,n}\to R_{2g+2}$ we obtain a universal curve on
$Y_{2g+2,n}$, with $n$
canonical sections.  The data of curve with sections defines a
map  $i_{g,n}:Y_{2g+2,n}\to\M_{g,n}$.
The universal curve on $\M_{g,n}$ (resp. $Y_{2g+2,n}$) gives rise to
a sheaf of relative first cohomology groups, and therefore to a continuous
$\ell$-adic representation $\pi_M$ of $\pialg(\M_{g,n},\bar m)$
(resp. $\pi_Y$ of $\pialg(Y_{2g+2,n},\bar y)$.
If we choose $\bar m$ to be the image
of a fixed geometric point $\bar y$ of $Y_{2g+2,n}$, we obtain the
commutative diagram

$$\tightmatrix{\pialg(Y_{2g+2,n},\bar y)&\rt{i_{g,n}}
&\pialg(\M_{g,n},\bar m)\cr
&{\scriptstyle\rho_Y}\rlap{$\searrow$}
\quad\swarrow\rlap{$\scriptstyle \rho_M$}&\cr
&\GL_{2g}(\Q_\ell)&\cr}$$

\ss The representation
$(\pi_Y,V_Y)$ is obtained by pull-back from the continuous representation
of $\pialg(R_n,\bar m)$ on the hyperelliptic curve $X_{\bar m}$.
Therefore, the geometric monodromy group $G_Y$ of the sheaf of
relative $H^1$ on $Y_{2g+2,n}$ is a subgroup of $\Sp(2g)$.
There are several ways to see that it is actually the full group.
One way is to construct a multi-section $R_n\to Y_{2g+2,n}$.
Another is to note that the inequality of exponential sums
$$(q+1-2g\sqrt{q}-n)^n\!\!\sum_{\a\in R_n(\F_q)}T_{\a}^{2k}\le
\sum_{(\a,x_1,\ldots,x_n)\in Y_{2g+2,n}(\F_q)}T_{\a}^{2k}\le
(q+1+2g\sqrt{q}-n)^n\!\!\sum_{\a\in R_n(\F_q)}T_{\a}^{2k}$$
implies
$$\sum_{(\a,x_1,\ldots,x_n)\in Y_{2g+2,n}(\F_q)}T_{\a}^{2k} = F(2k)
q^{n+2g+k} + O(q^{n+2g+k-1/2}),$$
and therefore
$$\dim\bigl({V_Y^{\t 2k}}^{G_Y}\bigr) = F(2k) =
\dim\bigl({V_Y^{\t 2k}}^{\Sp(2g)}\bigr).$$
It follows from a standard result in invariant theory \DMOS\ I Prop. 3.1 (c)
that $G_Y = \Sp(2g)$.

\proc{Theorem}{Let $p:X\to\M_{g,n}$ denote the universal
curve of genus $g$.  Then the geometric monodromy $G_{g,n}$
of $R^1 p_*\Q_\ell$ is $\Sp(2g)$.}

\proof The map $i_{g,n}$ realizes the geometric monodromy group
$G'_{2g+2,n}$ of $H^1$ of the universal curve on $Y_{2g+2,n}$ as a
subgroup of $G_{g,n}$. On the other hand,
by Poincar\'e duality, $G_{g,n}$ is contained
in $\Sp(2g)$.  The theorem follows.

\remark{Remark}We do not state a monodromy 
result about the moduli space $M_g$ itself because it does
not admit a universal curve.  There is an algebraic {\it stack} $\M_g$,
and presumably the natural language in which this theorem should be
framed is that of smooth $\ell$-adic sheaves on stacks.  The
choice of scheme language was dictated by the lack of adequate references
on the foundations of the theory of stacks.
\biblio
\end